\documentstyle[11pt]{article}
\newcommand{\cG}{{\cal G}}

\newcommand{\cP}{{\cal P}}
\newcommand{\cS}{{\cal S}}

\newtheorem{theorem}{Theorem}[section]
\newtheorem{lemma}{Lemma}[section]

\newtheorem{definition}{Definition}[section]
\newtheorem{corollary}{Corollary}[section]
\newcommand{\proof}{\noindent{\bf Proof.\ }}
\newcommand{\qed}{\hspace*{.25in} \rule{.75pt}{8pt}\hspace*{.03in}\rule{2pt}{8pt}}

\begin{document}
\title{Inverting Sets And The Packing Problem}

\author{Vance Faber\thanks{Work was performed under the auspices of the U.S. Department of Energy under Contract No. W-7405-ENG-36.} \\ 
Computer and Communication Division\\
Los Alamos National Laboratory\\
\\
Mark K. Goldberg\thanks{Work was
supported in part by the NSA under grant MDA904-90-H-4027 and by the NSF under grants IRI-8900511 and CDA-8805910; this work was performed while the author was visiting Los Alamos National Laboratory.}\\
Rensselaer Polytechnic Institute\\
Computer Science Department \\ 
\\
Emanuel Knill\\
School of Computer Science\\ 
Technical University of Nova Scotia\\
\\
Thomas H. Spencer\thanks{Work was supported in part by University Committee on Research, University of Nebraska at Omaha, and by the NSF under grants CDA-8805910 and CCR-8810609.}\\
Department of Mathematics and Computer Science \\
University of Nebraska at Omaha}
\date{}
\maketitle
\begin{abstract}

Given a set $V$, a subset $S$, and a permutation $\pi$ of $V$, we say
that  $\pi$ permutes $S$ if $\pi (S) \cap S = \emptyset$.
Given a collection $\cS = \{V; S_1,\ldots , S_m\}$, where
$S_i \subseteq V ~~(i=1,\ldots ,m)$, we say that $\cS$ is invertible if
there is a permutation $\pi$ of $V$ such that $\pi (S_i) \subseteq V-S_i$.
In this paper, we present necessary and sufficient conditions for the 
invertibility of  a collection and construct a polynomial algorithm 
which determines whether a given collection is invertible.  For an 
arbitrary collection, we give a lower bound for the maximum number
of sets that can be inverted. Finally, we consider the problem of
constructing a collection of sets such that no sub-collection
of size three is invertible. Our constructions of such collections
come from solutions to the packing problem with unbounded block sizes. We prove several new lower and upper bounds for the packing
problem and present a new explicit construction of packing.
\end{abstract}

\section{Introduction}

\par
The notion of invertibility arises as a tool for approaching
other combinatorial problems,  for example, the problem of
constructing a minimal size edge-set in a Cayley graph which intersects
every cycle of a given length.  We explain this connection using the hypercube as an example.

\subsection{Square-Blocking Edge-Sets in a Hypercube}

\par
Let $f(n)$ denote the minimal number of edges in a hypercube $Q_n$ such that
their removal from the hypercube yields a square-free graph. Evaluating 
$f(n)$ or computing the asymptotics for it is a long-standing open problem.
Erd\"{o}s (\cite{E4}) conjectured that $f(n)\asymp n2^{n-2}$ (see also
\cite{BS}, \cite{Ch}, \cite{D1}, \cite{D2}, \cite{GH}, \cite{JE}, \cite{HL}).
The best known
lower bound on $f(n)$ has been obtained by F.R.K. Chung~\cite{Ch}
and is given by $f(n)\geq (\alpha-o(1))n 2^{n-1}$, where
$\alpha$ is about $0.377$.

For every $n\geq 0$, we view $Q_{n+1}$ as the union of two 
copies of $Q_n$, denoted
respectively $Q'$ and $Q''$, and let $W$ denote the set of edges between
vertices in $Q'$ and vertices in $Q''$.  
Obviously, $W$ is a matching with $2^n$ edges; it can be viewed as a 
one-to-one mapping from $V(Q'')$ onto $V(Q')$
as well as  from $V(Q'')$ onto $V(Q')$. These mappings are naturally
expanded to one-to-one mappings of the corresponding edge-sets.  If
$K \subseteq E(Q'')$, then $W(K)$ denotes the image of $K$ 
under the mapping $W$. 
Given $W' \subseteq W$, let $C(W')$ denote  the set of vertices in 
$Q'$ that are incident to the edges in $W'$. 
Then the following statement can be easily proved.

\vskip 1\baselineskip
\noindent
{\bf Proposition.}
Let $N' \subset E(Q'),~N'' \subset E(Q'') $, and $W' \subset W$.
Then $N' \cup N'' \cup W'$ is a square-blocking set in $Q_{n+1}$ if and only if
$N'$ and $N''$ are square-blocking in $Q'$ and $Q''$, respectively, and
$C(W')$ is a vertex cover of the subgraph $E(Q') - N' - W(N'')$.
\qed

Since $Q_n$ is a connected bipartite graph, either of the partitions can 
serve as a vertex cover $C(W')$ for the corresponding $W' \subset W$.
With that, the choice of $N'$ and $N''$ is arbitrary, as long as both
are square-blocking sets for $Q'$ and $Q''$, respectively. Depending on
the choice of the construction for small values of $n$, this
construction leads to a square blocking set of size $(n-2)2^{n-2}$.

A smaller square-blocking set can be obtained if we try to construct 
$N'$ and $N''$ to minimize the intersection of $N' \cap W(N'')$. 
Thus, we may try to construct $N''$ to be an image of $N'$ under some permutation 
of the hypercube.  One way to construct such a permutation is to focus 
on the vertices of $Q'$ that are incident to at least $n/2$ edges in $N'$. 
Let $v' \in Q'$, $N'(v') \geq n/2$ and $v'' =W(v')$. Let 
$E'$ (respectively $E''$) be the edges in $N'$ (respectively in $N''$)
that  are adjacent to $v'$ (respectively $v''$).
We will be able to save at least one edge in $W'$ if 
$$F'\cap F' = \emptyset,$$
where $F' = E- E'$ and $F'' = E- E''$. In order to save more
edges, we may try to find a permutation of  all directions in $Q$
such that the condition above holds for as many vertices as possible.
If there are $m$ vertices $v_1,\ldots, v_m$ whose degrees $\geq n/2$,
and $S_i ~(i=1,\ldots,m)$ is the set of directions of $Q_n$ for which the
corresponding edges  are not in $N'$, then our goal is to permute
the set of all directions so that as many of $S_i's$ are inverted as
possible.

In Section 2, we establish a necessary and sufficient condition for 
a given collection of sets to be inverted by a single permutation,
and  give a bound for the number of sets that can be inverted in 
every collection $\cS$ with a given distribution of set sizes. 
An interesting reverse question is as follows: how many sets can a collection
have if no sub-collection containing a given number is invertible?
It turns out that even when the bound is three, the number is exponential, e.g.,
there are exponentially large collections of sets such that the maximal invertible
sub-collection contains  only two\footnote{It is easy to see that a collection 
of any two subsets of size not more than $n/2$ is invertible.} sets.
Construction of such a collection comes from solution of a packing problem 
for which the block sizes are unbounded.

\subsection{Packing Problem}

\par
Let $n$ be a positive integer, and let $c$ and $\alpha$ be two reals
in the interval $[0,1]$. The version of the packing problem
considered in this paper is defined as the construction of the  the maximum 
number of $cn$-subsets of a set
$[1,\ldots,n]$, such that any two of the sets intersect in fewer
than $\alpha c n$ elements.

Proving that a certain ``big'' packings exist is often done using the 
probabilistic method. The only algorithm that the method suggests is then 
the brute-force search, which is doubly-exponential for this problem.
For example, here is a result related to
Shannon's theorem (Z. F\"{u}redi, private communication):

\begin{theorem}
For every $0<d,c <1$, if $d>c^2$, then there exists $\epsilon>0$ such 
that there is a collection of $(1+\epsilon)^n$  $cn$-sets such that
the intersection of any two contains $<dn$ elements.
\end{theorem}

\proof 
({\it sketch}) Let $(1+\epsilon )^n = N$. Choose randomly 
$N$ sets of size $cn$ each. Then the expected size of the intersection of
any two is $c^2n$. By Chernoff's inequality 
$$Prob[ ~|A \cap A'| = c^2n +x \sqrt{n}] 
\stackrel{\displaystyle <}{\displaystyle \sim}
\exp (-x^2).$$
Since
$$Prob[ ~|A \cap A'| > dn ] \leq \exp(-(d-c)^2n),$$
one can delete a small number of bad elements.\qed

Since the probabilistic proof above does not provide an efficient way to
construct a packing, it is reasonable to seek algorithms that would
construct a packing with a sufficiently large number of blocks.
Since the output would have an exponentially large collection of sets,
the running time of the algorithm is inevitably an exponential function of
$n$.
Finally, if a solution to the packing problem is given by an explicit construction
we may expect that  the size of the collection is even smaller than  that
guaranteed by an algorithm.

\begin{definition}
A graph $\cG(n,c, \alpha)$ is defined as follows:
the vertices of $\cG$ are the $cn$-subsets  of $[1,\ldots,n]$;
two vertices are adjacent if and only if the corresponding sets intersect in
$\geq \alpha cn$ elements.  The number of vertices, the degree of a vertex, 
and the maximal size of an independent set of $\cG$ is denoted 
$N=N(n,c,\alpha),~ D =D(n,c,\alpha), $ and $P=P(n,c,\alpha)$, respectively.
\end{definition}
The Packing  Problem is then to evaluate the size of
the maximum independent set in $\cG(n,c,\alpha)$.

\begin{definition}
Given a set $V$, a collection $ \cS =S_1,\ldots, S_m$ of subsets of $\,V$,
and a set $\Gamma$ of permutations of $\,V$, $\kappa(\cS, \Gamma)$ is defined 
to be the maximum $k$ such that there exists a permutation $\pi \in \Gamma$ 
which inverts $k$ members of $\cS$.
If $\,\Gamma$ is the set of all permutations of $\,V$, then we write $\kappa(\cS)$
instead of $\kappa(\cS, \Gamma)$.
Given $S \subset V$ and a set $\Gamma$ of permutations of  $V$, 
$\lambda(S, \Gamma)$ denotes the number of permutations in $\Gamma$ 
that invert $S$.
\end{definition}

\section{Inverting Subsets of a Given Set}

\par
It turns out that there is a simple, necessary, and sufficient condition
for a collection $\cS =\{V; S_1,S_2,\ldots,S_m\}$ to be invertible.
We define a bipartite graph $G= G(\cS)$ with a bipartition $(V_1, V_2)$
as follows. Each of the sets $V_i~~(i=1,2)$ is in one-to-one
correspondence with $V$; two vertices $i\in V_1$ and $j\in V_2$ are
adjacent if and only if no set $S_k,~~(k=1,\ldots,m)$ contains both $i$ and $j$.

\begin{theorem}
A collection $\cS= \{V; S_1,\ldots, S_m\}$  is invertible if and only if
$G(\cS)$ has a perfect matching.
\end{theorem}
\proof 
If $\cS$ is invertible and $\pi$ is a permutation which
inverts each $S_i$'s, then  $\pi$ can also be viewed as the perfect
matching of $G$. The reverse is also straightforward.
\qed

\vskip 1\baselineskip\noindent
Two immediate  corollaries from the theorem above are:

\begin{corollary}
There is a polynomial algorithm which checks if a collection is
invertible, and if it is, outputs an inverting permutation. \qed
\end{corollary}

\begin{corollary}
If sets $\{S_1,\ldots, S_m\}$ are disjoint subsets of $V$, then
$\cS$ is invertible if and only if for every $i=1,\ldots,m,~|S_i|\leq |V|/2$. 
\end{corollary}
\proof 
Use the previous theorem together with the K\"{o}nig
condition on bipartite graphs with a perfect matching. \qed

\begin{corollary}
If $|V|=2k$, then $\cS = \{V,S_1,...,S_m\}$ with $|S_i|=k$ $(i=1,...,m)$ is invertible if and only if
$$|\cap_{i\in I} S_i\cap (\cap_{i\in I} \bar{S}_i)|=|\cap_{i\in I}\bar{S}_i\cap (\cap_{i\in I}S_i)|.$$
for every index set $I\subseteq[1,m]$.
\end{corollary}
\proof
If $\pi$ inverts $\cS$, then $\pi^{-1}$ inverts $\{ V,\bar{S}_1,...,\bar{S}_m\}$.  Thus \hfill\linebreak
$|\cap_{i\in I}S_i|=|\cap_{i\in I}\bar{S}_i|$ for any index set $I$.  The corollary follows by inclusion-exclusion.

\vskip 1\baselineskip
\noindent{\bf Remark.}  There are only at most $|V|$ non-empty conditions in Corollary 2.3.

\begin{theorem}
Let $\cS = \{V;S_1,S_2,S_3\}$ be an invertible collection satisfying
$|S_i| = k~(i=1,2,3)$. Then  $\cS$ is invertible if and only if
$$|S_1\cap S_2 \cap S_3| \leq 
|\overline{S_1} \cap \overline{S_2} \cap \overline{S_3}| \leq
|S_1\cap S_2 \cap S_3| + \frac{3}{2}(|V|-2k). \eqno{(*)}$$
\end{theorem}
\proof 
The proof uses inclusion-exclusion to check the conditions of Theorem 2.1.  We omit the details.

\begin{corollary}
Let $\cS$ be a collection of sets such that every
element of $V$ belongs to at most two sets of the collection. Then
$\cS$ is invertible if and only if every three sets in $\cS$ satisfy the condition of (\raisebox{-.04in}{*}) of Theorem 2.2.
\end{corollary}

\subsection{Inverting Large Sub-Collections}

\par
Our goal is to establish a bound on the number of sets in a given
collection $\cS$  that can be inverted by a single permutation. We
prove the existence of a permutation using a standard counting technique.
The  crucial detail in our case is that we consider a special class of
permutations, so called simple permutations. This restriction
substantially increases the lower bound. Given $\cS=\{V;S_1,\ldots, S_m\}$,
and  a permutation $\pi$ of $V$, we denote $\kappa(\cS, \pi))$ the
number of sets in the collection that are inverted by $\pi$. Then
$\kappa(\cS)= \max_{\pi} \kappa(\cS, \pi)$. If $\Pi$ is a given class
of permutations of $V$, $\lambda(S, \Pi)$ denotes the number of
permutations in $\Pi$ that invert the set $S$.
A permutation $\pi$ of a set with $n$ elements is called {\it simple},
if it has $\lfloor n/2 \rfloor$ disjoint cycles of length two;  $\sigma(n)$ 
denotes the number of simple permutations of a set with $n$ elements.

By extending the proof in [9], we get the following

\begin{lemma}
\[ \sigma(n)= \frac{n!}{ 2^{\lfloor n/2 \rfloor} \lfloor n/2 \rfloor !}. \]
\end{lemma}

Using Lemma 2.1 and a simple counting argument, we have the following lemma.

\begin{lemma}
Let $S$ be a subset of $V$. Then there are 
\[\frac{(n-i)!}{ 2^{\lfloor n/2 -i \rfloor} \lfloor n/2 -i \rfloor !}\]
simple permutations that invert $S$, where $i= |S| \leq n/2$.
\end{lemma}

\begin{theorem}
Let a  collection  $\cS$ contain $m_i$ sets of cardinality $i$ ($i =
1,\ldots , \lfloor
n/2 \rfloor )$. Then
\[  \kappa(\cS) \geq \frac{\lfloor n/2 \rfloor !}{n!}
		\sum_{i=1}^{\lfloor n/2\rfloor }
		\frac{(n-i) 2^i}{\lfloor n/2-i\rfloor !} m_i.
\]
\end{theorem}
\proof
Let $\Pi$ be the class of simple permutations of $V$. Obviously,
\[ \kappa(\cS) \geq \kappa(\cS, \Pi) \geq 
\frac{1}{\sigma (n)} \sum_{\pi \in \Pi} \kappa(\cS,\pi).\]
On the other hand, 
\[ \sum_{\pi \in \Pi} \kappa(\cS, \pi) =  \sum_{S\in \cS} \lambda(S,\Pi) .\]
Using the two previous lemmas we get the following:
\begin{eqnarray*}
\kappa(\cS) &\geq& \frac{2^{\lfloor n/2\rfloor } (\lfloor n/2\rfloor )!}{n!} 
\sum_{i=1}^{\lfloor n/2 \rfloor}\frac{(n-i)!}{2^{\lfloor n/2 \rfloor -i} (\lfloor n/2 \rfloor -i)!} m_i \\ 
&=& \frac{ \lfloor n/2 \rfloor !}{n!} \sum_{i=1}^{\lfloor n/2 \rfloor} \frac{2^i (n-i)!}{\lfloor n/2 -i\rfloor !} m_i. \qed
\end{eqnarray*}

\begin{corollary} 
If there exists $\epsilon >0$ such that for all
$i=1,\ldots \lfloor n/2 \rfloor$, $m_i \geq (1-\epsilon) {n \choose i}$, then 
\[ \kappa (\cS) \geq (1- \epsilon ) 3^{n/2}. \]
\end{corollary}
\proof
Using the previous theorem,
\begin{eqnarray*}
\kappa(\cS) &\geq& \frac{(\lfloor n/2 \rfloor )!}{n!} (1-\epsilon )
\sum_{i=1}^{\lfloor n/2 \rfloor} 
(\frac{2^i( n-i)!}{\lfloor n/2 -i\rfloor !} \times 
\frac{n!}{i! (n-i)!}) \\
&=& (1-\epsilon ) \sum_{i=1}^{\lfloor n/2 \rfloor} 
{\lfloor n/2 \rfloor \choose i}2^i = (1-\epsilon)3^{n/2}.\qed
\end{eqnarray*}

\begin{corollary}
There exists a sequence $\{M_n\}$ of square-blocking sets in the
hypercube $Q_n$, such that $|M_{n+1}|\leq 2|M_n|+2^n-f2^{n/3})$ where 
$f$ is a constant.
\end{corollary}
\proof
Using the previous results, one can construct square-blocking sets in a
$Q_n$ that  save $O(2^{n/3})$ edges of the hypercube.

\vskip 1\baselineskip
\noindent
{\bf Remark.} Erd\"{o}s asks:

\noindent
Given $n>0$, what is the largest $m$ such that there exists an invertible 
collection $\cS = \{ V;S_1,...,S_m\}$ with $|V|=n$?
  
\subsection{Set Collections with No Three Invertible}

\par
Given that $n$ is even, what is the largest number of subsets of $[n]$ of size 
$n/2$ such that no three are invertible? It turns out that such collections 
can be exponentially large.

\begin{lemma}
Let $k< n/2$, $K= [1,n/2 - k]$, and  let $\cP = \{R_i\}$  be
a collection of $k$-subsets  of $[1,n] -K$ such that the intersection 
of any two of them contain $<k/3$ elements. Then
for the collection $\cS = \{K \cup R_i\}$, no sub-collection of size
three is invertible.
\end{lemma}
\proof
If $S_i = K \cup R_i$ $(i=1,2,3)$ is a collection of three sets from
$\cS$, then
\[|\bigcap_{i=1}^3 S_i| = \frac{n}{2} -k +|\bigcap_{i=1}^3 R_i|,~~ and
~|\bigcap_{i=1}^3 \bar{S_i}| = \frac{n}{2} +k -|\bigcup_{i=1}^3 R_i|.\]
On the other hand, since 
$|\bigcup_{i=1}^3 R_i| = 3k - \sum_{i,j}|R_i \cap R_j| + 
|\bigcap_{i=1}^3 R_i| > 2k+ |\bigcap_{i=1}^3 R_i|$,
we see that the necessary condition of invertibility from Theorem 2.2 
does not hold for the collection ${S_1,S_2,S_3}$.
\qed

\vskip 1\baselineskip
\noindent
Thus, any packing $\Pi(n/2, c, 1/s)$ with exponentially many blocks
implies the existence  of an exponentially large collection of sets such
that no three of them are invertible.

\section{Packing with Unbounded Blocks}

\subsection{Lower Bounds for the Packing Problem}

\par
While there is a vast literature devoted to the  packing problem with
bounded block sizes (see \cite{DS} for references), there has
been relatively modest progress in the area of packing with unbounded block sizes. 
As noted in the introduction, the packing problem $\Pi(n,c,\alpha)$ is
equivalent to evaluating the maximal size of an independent set in
the the graph $\cG(n,c,\alpha)$.  
Our first bound follows from Tur\'{a}n's theorem.

\begin{theorem}
(Tur\'{a}n \cite{Tu}).
Every graph with $N$ vertices and average degree $D$ contains an 
independent set of size $\geq N/(D+1).$
\end{theorem}
Below, we use $N$ and $D$ to denote the vertex number and vertex degree
of the graph $\cG(n,c,\alpha)$; $S(c)$ denotes $c^{-c} (1-c)^{-1+c}$
for a given $c  (0< c<1)$. 
\begin{lemma}
There exists $A>0$ such that 
\[ {n \choose cn} \asymp \frac{A}{\sqrt n} S^n(c).\]
\end{lemma}
\proof
Use the Stirling formula.\qed

\begin{lemma}
If $\alpha >c$, then there exists $q >1$ such that for every $i > \alpha cn,$
\[\frac{{cn \choose i}{n-cn \choose cn-i}}{{cn \choose i-1}{n-cn
\choose cn-i+1}} >q.\]
\end{lemma}
\proof 
The following transformations are readily checked.
\[\frac{{cn \choose i}{n-cn \choose cn-i}}{{cn \choose i+1}
{n-cn \choose cn-i+1}}  = 
\frac{(i+1)!~(cn-i-1)!^2 ~(n-2cn +i+1)!} {i!~(cn-i)!^2 ~(n-2cn +i)!}\]
\[=\frac{(i+1)(n-2cn+i+1)}{(cn-i)^2} > \frac{i(n-2cn+i)}{(cn-i)^2} 
>\frac{\alpha (1-2c+\alpha c)}{c(1-\alpha)^2} >1.\]
The last inequality is equivalent to $\alpha >c$. \qed

\begin{lemma}
Let $\alpha > c$. Then there exists positive constants $A'~{\rm and}~A''$ 
such that 
\[N \asymp  \frac{A'}{\sqrt n} ~\frac {1}{c^{cn} (1-c)^{(1-c)n}};\]

\[D \asymp \frac{A''}{n} ~\frac{(1 -c)^{(1-c )n}}
{\alpha^{\alpha cn} (1-\alpha)^{2(1-\alpha) cn}
c^{c(1-\alpha )n} (1-2c + \alpha c)^{(1-2c + \alpha c)n}}.\]
\end{lemma}
\proof
The asymptotic for $N$ follows directly from Lemma 3.1.
>From the definition of the graph $\cG(n,c,\alpha)$,
$D= \sum_{i\geq \alpha cn}{cn \choose i}{n-cn \choose cn-i}.$
Then by Lemma 3.2, the first term of the summation is
the largest, and every other term is at least a constant smaller than the
previous. Thus, up to a constant,
\[D= {cn \choose \alpha cn}{n-cn \choose cn-\alpha cn}.\]
Using Lemma 3.1 again,  we have
\[
\begin{array}{rcl}
D&=& \frac{1}{n}~
\left(\frac{1} {\alpha^{\alpha }(1-\alpha)^{1-\alpha}}\right)^{cn}~
\left(\frac{1}{(\frac{(1-\alpha)c}{1-c})^{\frac{(1-\alpha)c}{1-c}}
 (1-\frac{(1-\alpha)c}{1-c})^{(1-\frac{(1-\alpha)c}{1-c})}}\right)^{1-c}\\
&=& \rule{0pt}{24pt}\frac{1}{n}~\frac{(1-c)^{(1-\alpha )cn} (1-c)^{(1-c)n -(1-\alpha)cn}}
      {\alpha^{\alpha cn} ~(1-\alpha)^{(1-\alpha
)cn}~((1-\alpha)c)^{(1-\alpha )cn}) 
      (1-2c + \alpha c)^{(1 -2c +\alpha c)n} }\\
&=& \rule{0pt}{24pt}\frac{1}{n}~\frac{(1-c)^{(1-c)n}}%
{\alpha^{\alpha cn} ~(1-\alpha)^{2(1-\alpha )cn}~ c^{(1-\alpha)cn } 
      ( 1-2c + \alpha c)^{(1 -2c +\alpha c)n} }\qed
\end{array}\]

\vskip 1\baselineskip
\begin{theorem}
Let 
$$T(n,c,\alpha) = ~\frac{\alpha^{\alpha cn}~ (1-\alpha)^{2(1-\alpha) cn}~ 
(1-2c + \alpha c)^{(1-2c + \alpha c)n}} {c^{\alpha c n}~ (1-c)^{2(1-c)n}}.
$$ 
Then, if $\alpha > c$, then  there is a packing $\Pi(n,c,\alpha)$ 
with at least $T(n,c,\alpha)$ blocks.
\qed
\end{theorem} 

The next theorem  shows how to compute the value of $c$ which maximizes
$T(n,c,\alpha)$ for a given $\alpha$.
\begin{theorem}
Given $\alpha >0$, the value of $c$ which maximizes $T(n,c,\alpha)$
is the positive root  of the following equation
$$
\alpha^\alpha (1-\alpha)^{2(1-\alpha)}(1-c)^2=
	 c^\alpha (1-2c+\alpha c)^{2-\alpha}.
$$
\end{theorem}
\proof
Let $f(c) = \log(T(n,c,\alpha))/n$.  We have
\begin{eqnarray*}
f(c) &=& \alpha c \log \alpha +2 (1-\alpha)c\log(1-\alpha)\\
&+&(1-2c+\alpha c)\log(1-2c+\alpha c)\\ 
&-& \alpha c \log(c) -2(1-c)\log(1-c)\\
&=& \alpha c (\log \alpha -2 \log(1-\alpha)+\log(1-2c+\alpha c))\\
&+&c(2\log(1-\alpha)\\
&-& 2\log(1-2c+\alpha c)+2\log(1-c))\\
&+&\log(1-2c+\alpha c) -2\log(1-c)-\alpha c \log(c). 
\end{eqnarray*}
Isolating the terms that are multiples of $c$, we get
\begin{eqnarray*}
f(c) &=& \mbox{}\!\!c (\alpha \log (\alpha/c) + 2(1-\alpha) \log(1-\alpha)
	+(\alpha-2)\log(1-2c+\alpha c)\\
&&+ 2\log(1-c)+\log(1-2c+\alpha c) -2\log(1-c).  
\end{eqnarray*}
Note that if $c=\alpha$ or as $c \rightarrow 0$ the logarithm is $0$.
Thus the maximum is in the range $0 < c < \alpha$.
After differentiation  and simplification,
\begin{eqnarray*}
f'(c)&=&
\left(\alpha \log (\alpha/c) + 2(1-\alpha) \log(1-\alpha)\right.\\
	&+& \left (\alpha-2)\log(1-2c+\alpha c)+2\log(1-c)\right) .
\end{eqnarray*}
We can rewrite this as
$$
f'(c)=\log\left(\alpha^\alpha (1-\alpha)^{2(1-\alpha)}(1-c)^2
	\over c^\alpha (1-2c+\alpha c)^{2-\alpha}\right).
$$
Thus $f'(c)=0$ when
$$
\alpha^\alpha (1-\alpha)^{2(1-\alpha)}(1-c)^2=
	 c^\alpha (1-2c+\alpha c)^{2-\alpha}. \qed 
$$ 

\begin{corollary}
For any particular $\alpha$ we can find the optimum $c$ by solving the 
equation above numerically. For example, if $\alpha= 1/3$, then
the optimal value of $c$ is close to $0.082508$, which yields the
base of the exponent in $T(n,c,\alpha)$ close to $1.0245$.
\end{corollary}

\subsection{Explicit Constructions}

\par
We would like to construct a family of sets without big intersections,
instead of just proving that such a thing exists.  

Recall that we are interested in the packing problem $\Pi(n,c, \alpha)$.
Let us assume that $k=1/\alpha$ is an integer.  One way to recursively
construct a packing is to divide the $n$ elements into $2k$ equal sized
disjoint subsets $A_i$.  Now we recursively construct packings 
$\Pi_i(n/2k, c, \alpha/2)$ on each subset $A_i$.  Each set $S_j$ of the
new packing is the union on one set from each of the $\Pi_i$.
We choose the $S_j$ so that no two of them have more than one set in
common.  The intersection of two of these sets has size at most
$$cn/(2k)+2k(n/2k)c(\alpha/2)=c\alpha n/2+nc\alpha/2=cn\alpha.$$

How many sets did we construct? First we need to count the number of $S_j$ 
as a function of $|\Pi_j|$.  Clearly the upper bound is $|\Pi_j|^2$.
In most cases this bound can be achieved.  Let $q=|\Pi_j|$ and assume that
there are more than $2k$ integers $a_i$ in the range $1 < a_i < q$ 
such that whenever $j \ne j'$ the difference $j-j'$ is relatively prime 
to $q$.  Number the elements of $\Pi_j$ 
from $0$ to $q$.  There will be one set $S_{lm}$ in the constructed family
$\Pi$ for each pair $(l,m)$ satisfying $0 \le l,m< q$.  The set 
$S_{lm}$ will contain the set numbered $l$ from $\Pi_1$ and the
set numbered $m$ from $\Pi_2$.  The contribution from $\Pi_j$ is
the set numbered $l+a_jm$ (provided that $j > 2$).

Now we want to see that two of the $S_{lm}$ share at most one set.
If they don't, then $l+a_jm=l'+a_jm'$ and $l+a_{j'}m=l'+a_{j'}m'$.
This is equivalent to $l-l'=a_j(m'-m)$ and $l-l'=a_{j'}(m'-m)$.
This means $a_j(m'-m)=a_{j'}(m'-m)$, so $m=m'$ since $a_j-a_{j'}$
is not a zero divisor.  Therefore, both sets are the same and no two
distinct $S_{lm}$ have two or more sets from the $\Pi_j$ in common.

Let $F(n,c,\alpha)$ represent the size of the
constructed family as well as the family itself.  As long as we avoid the
base case, we have
$$
F(n,c,\alpha)=(F(n\alpha/2,c,\alpha/2))^2.$$
We will say that the base case occurs when $n\alpha/4 \le 1$.
We will choose $c$ so that the base case construction consists of $n$
sets, each containing one element.  For some reason we are lead to conjecture
that the solution to this recurrence is 
$$\log(F(n,c,\alpha))=2^{A\sqrt{B\log n +C\log^2\alpha}+D\log \alpha}.$$
Substituting, we get
$$2^{A\sqrt{B\log n +C\log^2\alpha}+D\log \alpha}=
2^{A\sqrt{B\log(n\alpha/2) +C\log^2(\alpha/2)}+D\log(\alpha/2)+1}.$$
Taking logarithms and expanding, we get
\begin{eqnarray*}
&A&\!\!\!\!\!\!\!\!\sqrt{B \log n + C\log^2\alpha}+D\log \alpha\\ 
&\;\;\;=& \!\!A\sqrt{B\log n +B\log \alpha -B +C\log^2\alpha-2C\log \alpha +1}\\ 
&\;\;\;+& \!\!D\log\alpha -D+1
\end{eqnarray*}
If we choose $B=1$ and $C=1/2$, the right side simplifies to
\begin{eqnarray*}
&A&\!\!\!\!\!\!\!\!\sqrt{\log n + (1/2)\log^2\alpha}+D\log \alpha\\
&\;\;\;=& \!\!A\sqrt{\log n +(1/2)C\log^2\alpha}\\
&\;\;\;+& \!\!D\log\alpha -D+1
\end{eqnarray*}
Finally, we see that $D=1$ and that we can choose $A$ to make the base case
work.  Therefore
$$F=2^{2^{O(\sqrt{\log n})}}.\qed $$

\subsection{Upper Bounds}

\begin{theorem}
\label{TUB1}%
For $c>\alpha$, 
$P(n,c,\alpha)\leq {1-\alpha\over c-\alpha}$.
\end{theorem}

\proof
Let $\cS=\{X_1,\ldots, X_m\}$ be an independent subset of
$\cG(n,c,\alpha)$.
The size of $\bigcup\cS$ can be bounded from below by applying
the Schwarz inequality to indicator functions of sets
as done by Chung and Erd\"os in~\cite{CEe}:
\[
\Bigl(\sum_{1\leq i\leq m}|X_i|\Bigr)^2
  \leq |{\textstyle\bigcup}\cS|\sum_{1\leq i,j\leq m}|X_i\cap X_j|.
\]
This yields
\[
(mcn)^2\leq n(m(m-1)\alpha cn + mcn),
\]
and solving for $m$ gives $m\leq {1-\alpha\over c-\alpha}$.
\qed

The next result enables us to use Theorem~\ref{TUB1}
to obtain bounds on $P(n,c,\alpha)$ for $c\leq \alpha$.

\begin{theorem}
\label{TUB2}%
Let $1\geq e\geq c$ and $\alpha c\geq d\geq 0$. Then
\[
{(e-d)n \choose (e-c)n} P(n,c,\alpha)\leq {n \choose cn}
P\left((e-d)n,{c-d\over e-d}, {\alpha c-d\over c-d}\right).
\]
\end{theorem}

\proof
Let $\cS$ be an independent subset of $\cG(n,c,\alpha)$.
For $U\subset V\subseteq [1,\ldots, n]$ with
$|U|=dn$ and $|V|=en$, let
$\cS(U,V)=\{X\in \cS\;|\;U\subseteq X\subseteq V\}$.
For $X,Y\in \cS(U,V)$, we have
$|(X\setminus U)\cap(Y\setminus U)|< (\alpha c-d)n$.
Hence $|\cS(U,V)|\leq P((c-d)n,{c-d\over e-d},{\alpha c-d\over c-d})$.
There are ${n\choose en}{en\choose dn}$ many choices
for $U\subset V\subseteq [1,\ldots,n]$ with
$|U|=dn$ and $|V|=en$. Each $X\in \cS$ is a member of
${(1-c)n\choose (e-c)n}{cn\choose dn}$ many
$\cS(U,V)$. This gives
\[
|\cS|{(1-c)n\choose (e-c)n}{cn\choose dn}\leq
    {n\choose en}{en\choose dn}P\left((c-d)n,{c-d\over e-d},{\alpha c-d\over
       c-d}\right).
\]
The binomial identity 
${i\choose j}{i-j\choose k-j}={i\choose k}{k\choose j}$
for $i\geq k\geq j$
implies that
\[
{{(1-c)n\choose (e-c)n}{cn\choose dn}\Bigm/{(e-d)n\choose (e-c)n}}=
  {{n\choose en}{en\choose dn}\Bigm/{n\choose cn}}
.\]
The result follows.
\qed

Let $N(c,\alpha)={1-\alpha\over c-\alpha}$.
Combining Theorems~\ref{TUB1} and~\ref{TUB2}
yields the following corollary:

\begin{corollary}
\label{TUB3}%
If ${\alpha c-d\over c-d}<{c-d\over e-d}$, then
\[
P(n,c,\alpha)\leq {{n\choose cn}}
    N\left({c-d\over e-d},{\alpha c-d\over c-d}\right)
		\Bigm/{(e-d)n\choose (e-c)n}.
\]
\end{corollary}

We can now obtain good asymptotic
bounds on $\log P(n,c,\alpha)$.
Let $I(x)=-x\log(x)-(1-x)\log(1-x)$.
Corollary~\ref{TUB3} implies
$$
{\log P(n,c,\alpha)\over n}\leq
 I(c)-(e-d)I\left({ {e-c\over e-d}}\right)+o(1),\eqno{(*)}
$$
provided that ${\alpha c-d\over c-d}<{c-d\over e-d}$.
Let $B(e,d)=(e-d)I({e-c\over e-d})$.
To minimize the bound on $\log P(n,c,\alpha)$,
we find the maximum of $B(e,d)$ with the given constraints.
Note that the constraints are linear in $d$ and $e$, and
$B(e,d)$ is increasing in $e$ and decreasing in $d$.
By continuity, the bound of (*) holds for
${\alpha c-d\over c-d}={c-d\over e-d}$ and is minimized
when this identity holds.
Let $e'={(1-\alpha) c\over e-c}$
and $d'={(1-\alpha) c\over c-d}$. In terms of $e'$ and $d'$,
the constraint is $1-d'={e'\over d'+e'}$, which
implies that $e'+d'=1$. Additional constraints
on $e'$ and $d'$ are obtained
from the inequalities $1\geq e\geq c\geq \alpha c\geq d\geq 0$.
If $d=0$, then $d'=(1-\alpha)$. If $e=1$, then $d'={1-2c+c\alpha\over 1-c}$.
We now have
\[
B(e,d)={c(1-\alpha)\over d'(1-d')}I(d')
\]
to be maximized for $1-\alpha\leq d'\leq {1-2c+c\alpha\over 1-c}$.
The function $I(d')\over d'(1-d')$ is given by
$
 -{\log(d')\over (1-d')}-{\log(1-d')\over d'},
$
which is the sum of two convex functions on $(0,1)$.
(To see that f(x)=$-{\log(1-x)\over x}$ is convex,
write $f(x)=\sum_{i\geq 1}{x^{n-1}\over n}$.)
It follows that $B(e,d)$ is maximized on the boundary.
Thus our best asymptotic bounds on $\log P(n,c,\alpha)$
are obtained from  
\[
{\log P(n,c,\alpha)\over n}\leq I(c)-
{c(1-\alpha)\over d'(1-d')}I(d')+o(1)
\]
with $d'=1-\alpha$ or $d'={1-2c+c\alpha\over 1-c}$.
The value of $d'$ which yields the smaller bound depends on $c$ and $\alpha$.
To compare this
to the lower bounds obtained earlier,
consider $c=0.0825$ and $\alpha={1\over 3}$.
Then $\exp({\log P(n,c,\alpha)\over n})\leq 1.0655+o(1)$.
For $\alpha={1\over 3}$, the largest bound obtained for all $c$ occurs
for $c=0.1476$ and gives
$\exp({\log P(n,c,\alpha)\over n})\leq 1.0766+o(1)$.

\section*{Acknowledgements} 

\par
The authors are grateful to Jeffrey Dinitz, Paul 
Erd\"{o}s, Zoltan F\"{u}redi, and Earl S. Kramer for useful discussions.

\end{document}